\documentclass[a4paper]{article}

\usepackage{graphicx}%
\usepackage{multirow}%
\usepackage{amsmath,amssymb,amsfonts}%
\usepackage{amsthm}%
\usepackage{mathrsfs}%
\usepackage{xcolor}%
\usepackage{textcomp}%
\usepackage{manyfoot}%
\usepackage{booktabs}%
\usepackage{algorithm}%
\usepackage{algorithmicx}%
\usepackage{algpseudocode}%
\usepackage{listings}%

\usepackage{hyperref}
\hypersetup{unicode,
	colorlinks=true}
\usepackage[margin=1in]{geometry} 
\title{Discovering Power Grid Dynamics from Data Using Low-Rank Sparse Modeling}

\author{Aiman Mushtaq Purra\thanks{Department of Electrical Engineering, Institute of Technology, University of Kashmir, Srinagar, Jammu \& Kashmir, India 190006 (\url{aimanm1232@gmail.com},\url{danishrafiq@uok.edu.in})} \and Danish Rafiq\footnotemark[1]}

\begin{document}
\maketitle 

\begin{abstract}
The growing integration of renewable energy sources has significantly reduced grid inertia, making modern power systems more vulnerable to instabilities. Accurate estimation of dynamic parameters such as inertia constants and damping coefficients is critical, yet traditional model-based methods struggle with scalability and adaptiveness in large, low-inertia networks. This paper presents a novel data-driven framework that integrates Singular Value Decomposition (SVD) with Sparse Identification of Nonlinear Dynamics (SINDy) to estimate system parameters directly from time-series data. By reducing dimensionality before applying sparse regression, the proposed Latent-SINDy (L-SINDy) method mitigates overfitting while preserving essential system dynamics. The framework is validated on IEEE benchmark systems, including the 118-bus, 300-bus, and a large-scale 2869-bus European grid. The results demonstrate an accurate recovery of the inertia and damping values, with identified models that closely match the true dynamics of the system. The approach offers scalability, interpretability, and computational efficiency, highlighting its potential for real-time monitoring and control in renewable-rich grids.
\end{abstract}

\section{Introduction}
\subsection{Problem Statement}
The increasing penetration of renewable energy sources in modern power systems is fundamentally transforming grid dynamics \cite{enas2022inertia}. Conventional synchronous generators inherently provide system inertia through their rotating masses; however, these are increasingly being replaced by inverter-based resources that lack physical inertia \cite{abouyehia2025evaluating}. As a result, power systems are becoming more susceptible to frequency instabilities arising from reduced inertia \cite{candes2006robust}. Furthermore, key dynamic parameters such as the inertia constant and damping coefficient are no longer static but vary continuously with changing operating conditions \cite{ulbig2014impact}. The absence of accurate physical models for these evolving systems further complicates real-time modeling, analysis, and control \cite{borsche2015effects}. Consequently, dynamic parameter estimation has become essential for ensuring system stability. The core challenge lies in estimating these parameters adaptively without relying on fully known physical models, thereby motivating the use of data-driven approaches capable of learning system dynamics directly from measurements \cite{rudy2017data}.
\par
A considerable body of research has addressed the challenges associated with low-inertia power systems \cite{poudyal2021multiarea,oyewole2020power,hamid2022power,milano2018foundations}. Nevertheless, many existing techniques either presume prior knowledge of system equations or become computationally prohibitive for large-scale networks \cite{champion2019data}. Traditional state estimation methods are also inadequate when parameters vary frequently or remain uncertain \cite{luo2023unraveling}. Furthermore, most conventional modeling approaches rely heavily on detailed physical representations and fixed parameters, which often fail to capture real-time system behavior \cite{rosafalco2024ekf}. Although classical parameter estimation and machine learning-based methods have shown promise, they frequently suffer from overfitting, high computational complexity, and limited interpretability \cite{diggikar2025machine}. These limitations become more pronounced as system dimensionality increases. Recent studies have emphasized the importance of adaptive and scalable real-time models, particularly in systems with high renewable penetration \cite{rodales2023model}. However, most existing approaches fail to adequately capture nonlinear dynamics or scale efficiently to large interconnected networks \cite{raissi2019physics}. Hence, there is a pressing need for interpretable, data-efficient, and reduced-order modeling frameworks capable of handling the complexity of large-scale power systems.
\par
Accurate estimation of system parameters—specifically the inertia constants ($\mathbf{H}$) and damping coefficients ($\mathbf{D}$)—is crucial for stability assessment and control design \cite{zhao2019power}. In modern grids, these parameters evolve dynamically with operating conditions \cite{milano2018foundations}. Traditional model-based estimation approaches require detailed system equations and are computationally intensive, making them unsuitable for large networks \cite{zhao2019power}. The advent of Phasor Measurement Units (PMUs) provides high-resolution, time-synchronized data that offers new opportunities for data-driven modeling and parameter estimation \cite{bertozzi2024application, borsche2015effects}. Among various data-driven techniques, sparse regression methods—particularly the Sparse Identification of Nonlinear Dynamics (SINDy) framework \cite{brunton2016discovering,fasel2021sindy}—have shown great potential in identifying governing equations directly from measurement data. However, direct application of SINDy to large-scale power systems often results in overfitting and numerical ill-conditioning \cite{anguluri2022parameter}. This limitation underscores the need to incorporate dimensionality reduction before applying SINDy. Despite its potential, such an integrated framework has not yet been systematically explored for large-scale power grids.
\par
To address this gap, this paper proposes a hybrid modeling framework that combines dimensionality reduction with sparse identification for dynamic parameter estimation in large-scale, low-inertia power systems. The principal contributions of this work are summarized as follows:
\begin{enumerate}
\item We propose a hybrid framework that integrates Singular Value Decomposition (SVD) with the SINDy algorithm for dynamic parameter estimation in low-inertia grids.
\item The proposed method is applied to large IEEE benchmark systems, including the 118, 300, and 2869-bus test cases.
\item We demonstrate that the Latent-SINDy (L-SINDy) framework achieves accurate estimation of inertia constants and damping coefficients while mitigating overfitting in high-dimensional systems.
\item We provide an open-source implementation code, in the form of MATLAB scripts, to reproduce the results of the paper.
\end{enumerate}
\subsection{Power system dynamical model}
The swing equation is a fundamental model for analyzing the transient dynamics of synchronous generators in power systems. It describes how a generator’s rotor angle and speed evolve in response to disturbances such as load variations, faults, or generator outages that disrupt the balance between mechanical input power and electrical output power. Under nominal operating conditions, the mechanical torque from the prime mover equals the electrical torque delivered to the grid, maintaining a constant rotor speed. When this balance is disturbed, the rotor accelerates or decelerates depending on whether the mechanical power exceeds or falls below the electrical power. The swing equation captures this dynamic behavior by incorporating the effects of inertia, damping, and power imbalance during system disturbances.
\par
Mathematically, the swing equation is expressed as the nonlinear second-order differential equation:
\begin{equation}\label{swing}
\frac{2\mathbf{H}_i}{\omega_R} \ddot{\boldsymbol{\delta}_i} + \frac{\mathbf{D}_i}{\omega_R} \dot{\boldsymbol{\delta}_i} = \mathbf{F}_i - \sum_{j=1;j \neq i}^{n} \mathbf{K}_{ij} \sin(\boldsymbol{\delta}_i - \boldsymbol{\delta}_j - \boldsymbol{\gamma}_{ij})
\end{equation}
where $\boldsymbol{\delta}_i$ is the rotor angle of the $i$th generator, $\omega_R$ denotes the system reference angular frequency, and $\mathbf{H}_i$ and $\mathbf{D}_i$ represent the inertia and damping constants, respectively. The term $\mathbf{K}_{ij} \geq 0$ denotes the coupling between oscillators $i$ and $j$, and $\boldsymbol{\gamma}_{ij}$ is the corresponding phase shift in this coupling. The constants $\mathbf{F}_i$, $\mathbf{K}_{ij}$, and $\boldsymbol{\gamma}_{ij}$ are derived from power flow solutions and Kron reduction \cite{nishikawa2015comparative}.
\par
Three principal modeling approaches are widely used in power system studies: the Effective Network (EN), Synchronous Motor (SM), and Structure Preserving (SP) models. These models primarily differ in their representation of loads, which affects the interpretation of $\mathbf{F}_i$, $\mathbf{K}_{ij}$, and $\boldsymbol{\gamma}_{ij}$. In the EN model, loads are treated as constant impedances. The SP model represents loads as first-order oscillators with $\mathbf{H}_i = 0$, while the SM model treats them as synchronous motors represented by second-order oscillators. The inertia constant $\mathbf{H}$ plays a vital role in determining the rate at which a rotor responds to power imbalances—higher inertia values dampen speed deviations and enhance system stability. However, as renewable energy integration increases, system inertia naturally decreases, making the swing equation even more critical for assessing rotor stability and maintaining overall grid reliability in modern low-inertia systems.
\par
The remainder of this paper is organized as follows. Section \ref{sec-2} reviews related work on inertia reduction, model reduction, and data-driven parameter estimation. Section \ref{sec-LSINDy)} presents the proposed Latent-SINDy (L-SINDy) framework. Section \ref{sec-results} details the case studies conducted on IEEE benchmark systems and discusses the corresponding results. Finally, Section \ref{sec-conclusion} concludes the paper and outlines key limitations and directions for future research.

\section{Related Works}\label{sec-2}
A variety of methods have been developed for inertia estimation and dynamic model identification in power systems \cite{oyewole2020power,hamid2022power,milano2018foundations}. Classical approaches rely on analytical formulations, such as the Rate of Change of Frequency (RoCoF) or energy-based methods, and on \textit{Kalman} filtering techniques that assume a known model structure \cite{rosafalco2024ekf}. These techniques often require prior knowledge of system parameters and are limited in their adaptability to rapidly evolving grid conditions. Dynamic State Estimation (DSE), as reviewed in \cite{zhao2019power}, focuses on data-driven state tracking but still depends on parameter knowledge or full state observability.
\par
Recent advancements in data-driven methods have introduced Physics-Informed Neural Networks (PINNs) and other Machine Learning (ML) approaches for estimating generator parameters from measurement data \cite{raissi2019physics,diggikar2025machine}. However, such methods may exhibit slow convergence and require large datasets, particularly in low-inertia systems. On the model reduction front, techniques such as Proper Orthogonal Decomposition (POD) and Singular Value Decomposition (SVD) have been widely employed for large-scale systems. For instance, \cite{khatibi} compared modal, balanced, and SVD-based reduction techniques and demonstrated that reduced-order models can capture dominant electromechanical modes with significantly fewer states. Nevertheless, linear reduction methods may fail to represent nonlinear interactions in higher-order dynamics.
\par
Related research has also focused on data-driven equivalents and dynamic aggregation. The IEEE Task Force on DSE \cite{zhao2019power} has emphasized the importance of reduced-order models for enabling scalable state estimation. Despite these developments, few studies combine model reduction with sparse regression. Recent reviews further highlight that inertia estimation in low-inertia grids remains challenging and advocate hybrid approaches that integrate physics-based and data-driven techniques \cite{abouyehia2025evaluating,poudyal2021multiarea}.
\par
The Sparse Identification of Nonlinear Dynamics (SINDy) algorithm \cite{brunton2016discovering} performs sparse regression on a library of candidate nonlinear functions to infer governing equations directly from data. Compared to PINNs and Extended Kalman Filters (EKFs), SINDy offers lower computational complexity and better generalization to unseen conditions. It has been successfully applied to a broad range of systems, including fluid dynamics \cite{loiseau2018constrained}, optical devices \cite{sorokina2016sparse}, chemical reaction networks \cite{hoffmann2019reactive}, plasma convection \cite{dam2017sparse}, and structural systems \cite{lai2019sparse}, as well as in designing model predictive controllers \cite{kaiser2018sparse}.
\par
In addition to practical applications, several theoretical extensions of SINDy have been proposed, including formulations for identifying partial differential equation (PDE) systems \cite{schaeffer2017learning,rudy2017data}, rational nonlinearities \cite{mangan2017inferring}, incomplete datasets \cite{schaeffer2017sparse}, and models incorporating physical constraints \cite{loiseau2018constrained}. Moreover, SINDy has been combined with deep autoencoder architectures to discover intrinsic low-dimensional coordinates through joint optimization \cite{champion2019data}. However, its applicability to large-scale power systems has not yet been extensively explored. The following subsection reviews the SINDy methodology, while its implementation for dynamic modeling in large-scale power networks is presented in the next section.
\subsection{SINDy method}
Consider a nonlinear dynamical system represented in state-space form as:
\begin{equation}\label{FOM}
     \mathbf{\dot{x}}(t) = \mathbf{f}(\mathbf{x}(t)),
\end{equation}  
where $\mathbf{x}(t) \in \mathbb{R}^n$ is the state vector with $n$ number of state variables, and the function $\mathbf{f}(): \mathbb{R}^{n} \mapsto \mathbb{R}^{n}$ denotes the nonlinear vector field governing the state evolution. Let $m$ be the number of temporal snapshots and $t$ be the time vector $t = t_1, t_2, \ldots, t_m$. To approximate $\mathbf{f}$ from data, we construct the state matrix $\mathbf{X} \in \mathbb{R}^{n \times m}$ as:
\begin{equation}\label{Xmatrix}
	\mathbf{X} = \begin{bmatrix}
\mathbf{x}_1(t_1) & \mathbf{x}_1(t_2) & \cdots & \mathbf{x}_1(t_m) \\
	\mathbf{x}_2(t_1) & \mathbf{x}_2(t_2) & \cdots &\mathbf{x}_2(t_m) \\
	\vdots & \vdots & \ddots & \vdots \\
	\mathbf{x}_n(t_1) & \mathbf{x}_n(t_2) & \cdots & \mathbf{x}_n(t_m)	
	\end{bmatrix}.
\end{equation}
The corresponding time-derivative matrix $\mathbf{\dot{X}}$ is given by:
\begin{equation}\label{Xdot}
	\mathbf{\dot{X}} = \begin{bmatrix}
	\mathbf{\dot{x}}_1(t_1) & \mathbf{\dot{x}}_1(t_2) & \cdots & \mathbf{\dot{x}}_1(t_m) \\
	\mathbf{\dot{x}}_2(t_1) & \mathbf{\dot{x}}_2(t_2) & \cdots & \mathbf{\dot{x}}_2(t_m) \\
	\vdots & \vdots & \ddots & \vdots \\
	\mathbf{\dot{x}}_n(t_1) & \mathbf{\dot{x}}_n(t_2) & \cdots & \mathbf{\dot{x}}_n(t_m)	
	\end{bmatrix}.
\end{equation}
We define the library of nonlinear candidate terms as follows:
	\begin{equation}\label{theta}
\mathbf{\Theta(X)} = \begin{bmatrix}
	| & | &  | &  & | & | &  \\
	1 & \mathbf{X} & \mathbf{X}^{P_2}  & \cdots &  \sin\mathbf{X} & \cos\mathbf{X} & \cdots \\	
	| & | &  | &  & | & | &  \\
		\end{bmatrix},
	\end{equation}
where $\mathbf{X}^{P_2}$ represent higher-order polynomial terms. For instance, the quadratic polynomial candidates are given by:
\begin{align}\label{P}
	&\mathbf{X}^{P_2} = \nonumber\\
    &\begin{bmatrix}
	\mathbf{x}_1^{2}(t_1) &  \mathbf{x}_1(t_1)\mathbf{x}_2(t_2) & \cdots & \mathbf{x}_2^{2}(t_1) & \cdots & \mathbf{x}_n^{2}(t_1) \\
		\mathbf{x}_1^{2}(t_2) &  \mathbf{x}_1(t_2)\mathbf{x}_2(t_2) & \cdots & \mathbf{x}_2^{2}(t_2) & \cdots & \mathbf{x}_n^{2}(t_2) \\
		\vdots  & \vdots & \ddots &  \vdots & \ddots & \vdots \\
			\mathbf{x}_1^{2}(t_m) &  \mathbf{x}_1(t_m)\mathbf{x}_2(t_m) & \cdots & \mathbf{x}_2^{2}(t_m) & \cdots & \mathbf{x}_n^{2}(t_m) 
	\end{bmatrix}.	
\end{align}
SINDy then represented the nonlinear system (\ref{FOM}) in terms of the data matrix (\ref{Xmatrix}) and (\ref{Xdot}) as:
\begin{equation}\label{X}
   \mathbf{\dot{X}} = \mathbf{\Theta(X)}\mathbf{\Xi} ,
\end{equation}
where \( \mathbf{\Xi} \in \mathbb{R}^{p \times n} \) is a sparse coefficient matrix with each column encoding the active functions governing the dynamics of the corresponding state variable.
The $k$th state equation can thus be written as:
    \begin{equation}\label{eqn-SINDyk}
         \mathbf{\dot{x}}_k = \mathbf{f}_k\mathbf{(x)} = \mathbf{\Theta}\mathbf{(x)}^T\mathbf{\zeta}_k
    \end{equation}
where $\mathbf{\zeta}_k$ denotes the $k$th column of $\mathbf{\Xi}$. Consequently, the overall system can be reconstructed as:
\begin{equation}\label{system}
  \mathbf{\dot{x}}(t) = \mathbf{f}\mathbf{(x}(t)) = \mathbf{\Xi}^T(\mathbf{\Theta}\mathbf{x}(t)^T),
\end{equation}
indicating that system trajectories are represented by a sparse linear combination of nonlinear basis functions. The active coefficients in $\mathbf{\Xi}$ are identified by solving the following sparse regression problem:
\begin{equation}\label{sparse}	
\mathbf{\Xi} =		\min_{\Theta} \left\| \mathbf{\Theta(X)\Xi} - \mathbf{\dot{X}} \right\|_2^2 + \lambda \|\mathbf{\Xi}\|_1
\end{equation}
where the first term enforces data fidelity and the second term ($\ell_1$-regularization) promotes sparsity, ensuring that only the dominant terms are retained. This leads to a parsimonious yet accurate representation of the underlying system dynamics.
\subsection{Issues in SINDy}
The SINDy algorithm requires constructing the data matrix $\mathbf{\Theta}(\mathbf{X})$, containing candidate functions (monomials, trigonometric, etc.) of the state variable $\mathbf{x}\in\mathbb{R}^n$. In power systems, the snapshot matrix $\mathbf{X}$ is generally large due to the involvement of thousands of state variables, resulting in high dimensionality ($n>>1$). The large state dimension leads to factorial growth of $\mathbf{\Theta}(\mathbf{X})$, leading to an overdetermined solution of (\ref{X}) \cite{brunton2016discovering}. Even with scaling, a large $\Theta$ matrix imposes a heavy computational bottleneck and requires far more data samples to achieve a well-determined sparse regression. Additionally, each row of $n$  requires a separate optimization problem to be solved (cf. eq. (\ref{eqn-SINDyk})), which immediately leads to a computationally expensive task. Thus, applying SINDy directly to power grids can lead to conditioning issues, spurious nonlinear terms that do not represent known connectivity patterns, or complete failure. Therefore, there is a strong motivation to first obtain a low-rank approximation of the original state vector $\mathbf{x}$ using dimensionality reduction, i.e.,  $\mathbf{x}\approx\mathbf{\Phi}\mathbf{z}$, where $\mathbf{z} \in \mathbb{R}^{r}, r<<n$,
and then obtain a sparse representation in terms of coefficients of $\mathbf{z}$. This is explained in the next section.

\section{Latent-SINDy for power systems}\label{sec-LSINDy)}
\begin{figure}[H]
\centering
\includegraphics[width=\linewidth]{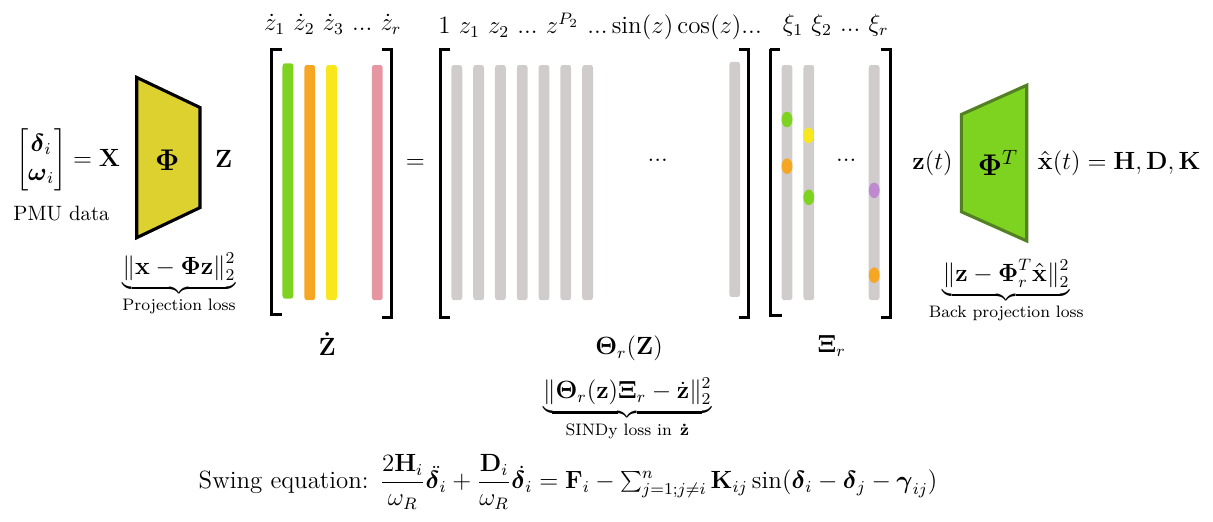}
\caption{Graphical illustration of LSINDy}
\label{fig:block diagram}
\end{figure}

In this section, we will describe how the SINDy method can be applied to the latent variables for power system applications. To obtain the nonlinear function $(\mathbf{f}$, in eq. (\ref{FOM})), from the power system data, we collect the state data matrix $ \mathbf{X} $ given as:
\begin{equation}\label{data}
       \mathbf{X} = \begin{bmatrix}
		\boldsymbol{\delta}(t_0) & \boldsymbol{\delta}(t_1) & \cdots & \boldsymbol{\delta}(t_m) \\
		\boldsymbol{\omega}(t_0) & \boldsymbol{\omega}(t_1) & \cdots &\boldsymbol{\omega}(t_m)	
	\end{bmatrix}   \in \mathbb{R}^{2n_g \times t_m}
\end{equation}
and \begin{equation}\label{derivative}
\mathbf{\dot{X}} = \begin{bmatrix}
		\boldsymbol{\dot{\delta}}(t_0) & \boldsymbol{\dot{\delta}}(t_1) & \cdots & \boldsymbol{\dot{\delta}}(t_m) \\
			\boldsymbol{\dot{\omega}}(t_0) & \boldsymbol{\dot{\omega}}(t_1) & \cdots & \boldsymbol{\dot{\omega}}(t_m) 
	\end{bmatrix} \in \mathbb{R}^{2n_g \times t_m},
\end{equation}
where $m$ is the number of snapshots collected, $n_g$ is the number of oscillators in the system and $\boldsymbol{\omega}$ is the angular frequency of each rotor. By the application of SVD, we obtain the following:
\begin{equation}\label{SVD}
    \mathbf{X} = \mathbf{U}\mathbf{\Sigma}\mathbf{V}^T 
\end{equation}
where $\mathbf{U} \in \mathbb{R}^{n \times n}$ and $\mathbf{V} \in \mathbb{R}^{t_m \times t_m}$ are orthogonal matrices and $\mathbf{\Sigma} \in \mathbb{R}^{n \times t_m}$ is a matrix with real, non-negative entries on the diagonal and zeros off the diagonal. The entries of $\mathbf{\Sigma}$, i.e., $\sigma_1, \sigma_2,...,\sigma_n$ are referred to as singular values and are ordered hierarchically from most significant to least significant, i.e., $\sigma_1> \sigma_2>...>\sigma_n$.
Thus, by truncating the $r$ columns of $\mathbf{U}, \mathbf{V}$ and $\mathbf{\Sigma}$ given as:
\begin{equation}\label{truncatedSVD}
\mathbf{X}=\mathbf{U\Sigma V^T}= \begin{bmatrix}
    \mathbf{{U}}_r& \mathbf{\hat{U}}^{\bot}
    \end{bmatrix}
    \begin{bmatrix}
        \mathbf{{\Sigma}}_r\\\mathbf{0}
    \end{bmatrix}\mathbf{V}^T=\mathbf{{U}}_r\mathbf{{\Sigma}}_r\mathbf{V}^T
\end{equation}
one can obtain a low-rank approximation of $\mathbf{X}$ if $r<n$, where the columns of $\mathbf{\hat{U}}^{\bot}$
span a vector space that is complementary and orthogonal to that spanned by $\mathbf{{U}}_r$. This yields the projection matrix $\mathbf{\Phi}_r \in \mathbb{R}^{n \times r}=\mathbf{U}_r$ by retaining first $r$ columns of $\mathbf{U}$. Using the projection matrix $\mathbf{\Phi}$, we can obtain the low-order approximation of the state vector $\mathbf{x}(t)$ as: \begin{equation}\label{approximation}
\mathbf{x}(t) \approx \mathbf{\Phi}_{r}\mathbf{z}(t),
\end{equation}
where $ \mathbf{z}(t) \in \mathbb{R}^{r} $ is the latent state variable. Using the approximation
(\ref{approximation}) in (\ref{FOM}), we can derive the reduced dynamical model, which we refer to as the reduced order model (ROM) of (\ref{FOM}), given as:
\begin{equation}\label{ROM}
\mathbf{\dot{z}}(t) = \mathbf{\Phi}_{r}^{T}\mathbf{f}_r(\mathbf{\Phi}_{r}\mathbf{z}(t)),
\end{equation}
where $\mathbf{f}_r:\mathbb{R}^n \mapsto \mathbb{R}^{r}$ is the reduced nonlinear function. It should be noted that since $\mathbf{f}$ and $\mathbf{f}_r$ remain unknown, we instead obtain the reduced snapshot matrix $\mathbf{Z}$ as:
\begin{equation}\label{reducedstates}
    \mathbf{Z}=\mathbf{\Phi}_{r}^{T}\mathbf{X} \in \mathbb{R}^{r \times t_m},
\end{equation} which is equivalent to expressing each full state-variable $\mathbf{x}(t_i)$ in the reduced orthonormal basis $\mathbf{\Phi}_r$. Thus, we have:
\begin{equation}\label{Z}
 \mathbf{Z} = \begin{bmatrix}
		 | & | &  & | &    \\
		\mathbf{z}(t_1) & \mathbf{z}(t_2) & \cdots & \mathbf{z}(t_m) \\
		 | & | &  & | &    \\
	\end{bmatrix},
\end{equation} and its derivative matrix 
\begin{equation}\label{Zdot}
	\mathbf{\dot{Z}} 
	= \begin{bmatrix}
			| & | &  & | &    \\
		\mathbf{\dot{z}}(t_1) & \mathbf{\dot{z}}(t_2) & \cdots & \mathbf{\dot{z}}(t_m) \\
			| & | &  & | &   \\
	\end{bmatrix}
\end{equation}
obtained as 
\begin{equation}\label{reducedstates}
\mathbf{\dot{Z}}=\mathbf{\Phi}_{r}^{T}\mathbf{\dot{X}} \in \mathbb{R}^{r \times t_m}.
\end{equation}
Then we perform the SINDy method on $ \mathbf{Z} $ and $ \mathbf{\dot{Z}}$ to obtain the latent dynamics model as:
\begin{equation}\label{ROMz}
\mathbf{\dot{z}}_k=\mathbf{\Theta_r}(\mathbf{z}^T)\boldsymbol{\beta}_k ,
\end{equation}
where $\boldsymbol{\beta}_k$ is the $k$th column of the reduced coefficient matrix $\mathbf{\Xi}_r \in \mathbb{R}^{q \times r}$, $q$ is the number of library functions used in the reduced coordinates $\mathbf{z}$, and 
\begin{equation}\label{theta_r}
		\mathbf{\Theta_{r,i}} = \begin{bmatrix}
			| & | & | & | &  & | & | &  \\
			1 & \mathbf{Z} & \mathbf{Z}^{P_2} & \mathbf{Z}^{P_3} & \cdots &  \sin\mathbf{Z} & \cos\mathbf{Z} & \cdots \\	
			| & | & | & | &  & | & | &  \\
		\end{bmatrix}.
	\end{equation}
To identify active terms in the latent space, we solve a sparse regression problem in $\mathbf{\Xi}_r$, given as:
\begin{equation}\label{Zsparse}
	\mathbf{\Xi}_r =		\min_{\Theta_r} \left\| \mathbf{\Theta_r(Z)\Xi}_r - \mathbf{\dot{Z}} \right\|_2^2 + \lambda \|\mathbf{\Xi}_r\|_1  \\	
\end{equation}
Finally, we can reconstruct the full-state prediction ${\mathbf{\hat{x}}}\approx\mathbf{\Phi}_{r}\mathbf{z}$ by integrating the identified reduced-order model (\ref{ROMz}).	
The above discussion is illustrated in Fig. \ref{fig:block diagram} and the steps are summarized in Algorithm (\ref{alg:L-SINDy}). In practice, we tune the SINDy regularization parameter $\lambda$ to balance sparsity and accuracy. The latent dimension $r$ can be chosen by inspecting the singular value decay or by cross-validation. In the next section, we will show the application of L-SINDy for power systems.
\begin{algorithm}[H]
\caption{Latent Sparse Identification of Nonlinear Dynamics (L-SINDy)}\label{alg:L-SINDy}
\begin{algorithmic}[1]
\Statex \textbf{Input:} Rotor-angle $\delta_i$ and frequency data $\omega_i$
\Statex \textbf{Output:} Sparse coefficient matrix $\mathbf{\Xi_r}$, reconstructed state vector $\mathbf{\hat{x}}$
\State Formulate $\mathbf{X}$ as per eqn. \eqref{data}
\State [$\mathbf{U},\mathbf{\Sigma},\mathbf{V}]=svd(\mathbf{X})$ \Comment{perform SVD on $\mathbf{X}$}
\State $\mathbf{\Phi_r}=truncate(\mathbf{U})$\Comment{retain $r$ columns of $\mathbf{U}$}
\State Formulate $\mathbf{Z}$ and $\mathbf{\dot{Z}}$ as per  \eqref{Z} and (\ref{Zdot})
\State Construct the library $\mathbf{\Theta_r}_i$ as per eqn. \eqref{theta_r}
\State Set tuning parameter $\lambda$
\State Perform sparse regression  \eqref{Zsparse} to obtain $\mathbf{\Xi_r}$
\State Solve the reduced system $\mathbf{\dot{z}}=\mathbf{\Theta_r}(\mathbf{z}^T)\boldsymbol{\beta}_k$ 
\State Back project to obtain the full state dimension ${\mathbf{\hat{x}}}\approx\mathbf{\Phi}_{r}\mathbf{z}$
\end{algorithmic}
\end{algorithm}
\section{Numerical Simulations}\label{sec-results}
We apply the proposed L-SINDy method to three test systems of increasing size: (i) \textit{IEEE 118-bus} having 54 generators, (ii) \textit{IEEE 300-bus} with 69 generators, and (iii) \textit{2869-bus European high-voltage} network having 510 generators. The simulations have been performed in MATLAB version R2018a running on an 11th Gen Intel(R) Core(TM) i5-11300H @ 3.10GHz desktop computer with 32 GB RAM. All models have been simulated using MATLAB's built-in \textit{ode45} solver \cite{DormandPrince1980} for a time span of $[0,5s]$ with a step size of $\Delta t =0.01s$. The data was obtained from the MATPOWER library \cite{zimmerman2010matpower} to construct the EN power grid models as described in \cite{nishikawa2015comparative}.
\subsection{IEEE-118 Bus system}
\begin{figure}[htbp]
\centering
\includegraphics[width=0.5\linewidth]{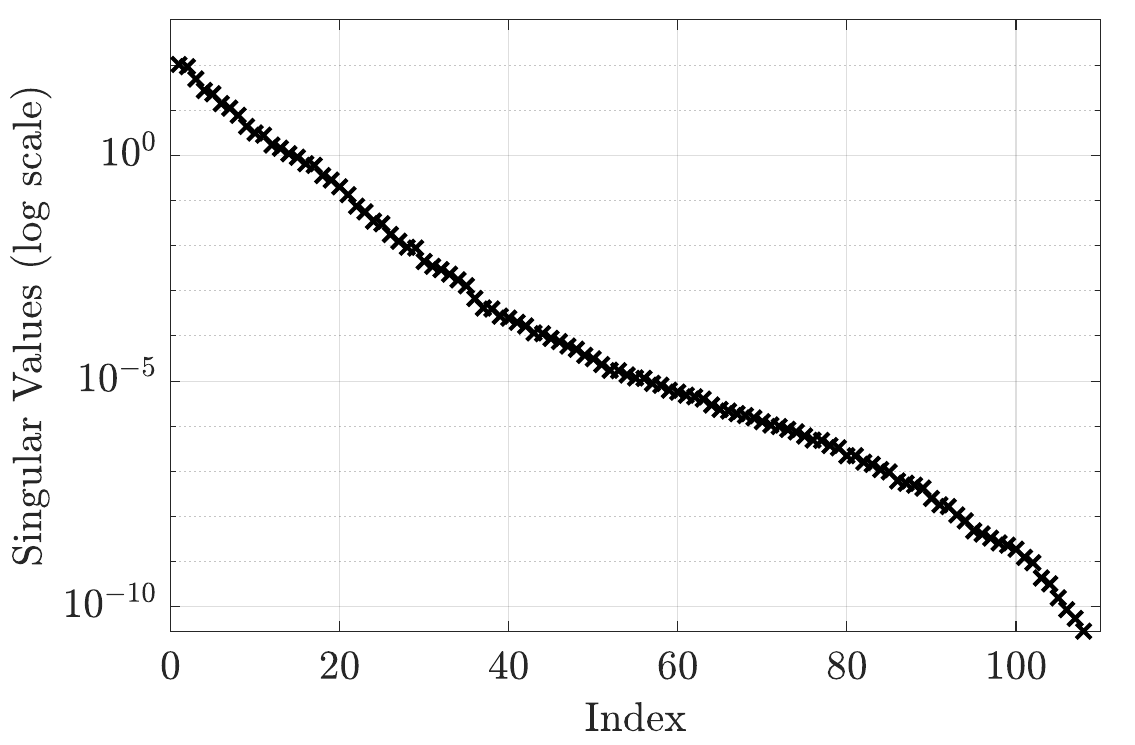}
	\caption{Singular Value decay plot of IEEE-118 bus system}
	\label{fig:decay118}  
\includegraphics[width=0.5\linewidth]{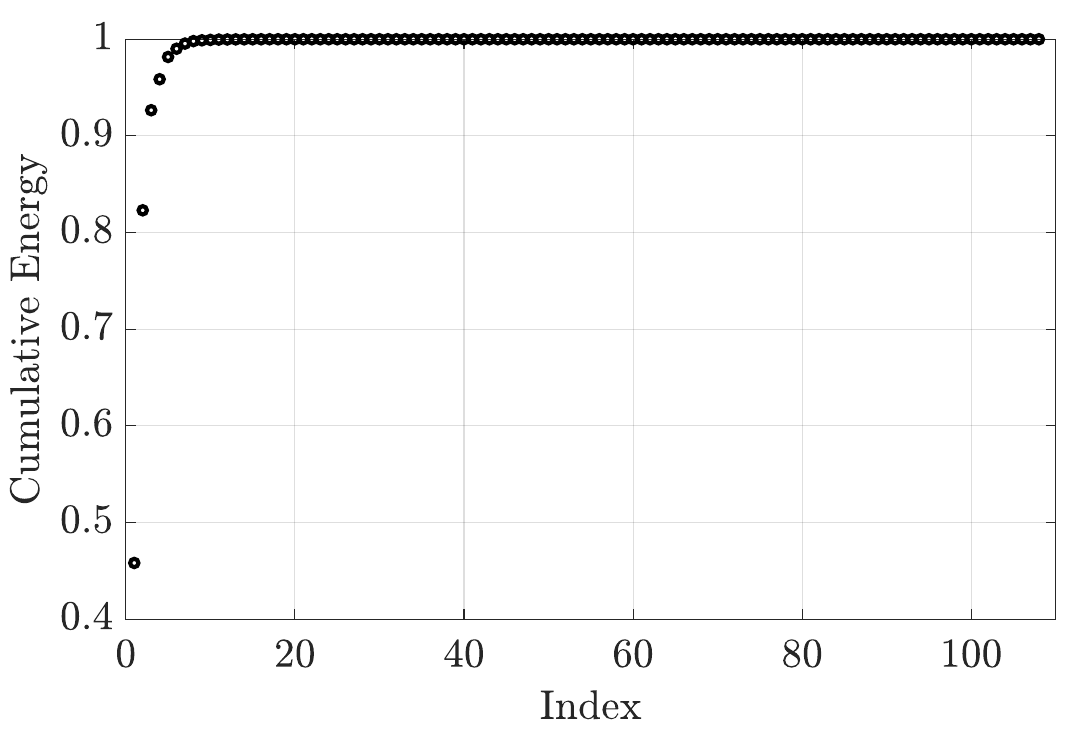}
	\caption{Singular Value energy plot of IEEE-118 bus system}
	\label{fig:energy118}	
\end{figure}

\begin{figure}[H]
\centering
\includegraphics[width=\linewidth]{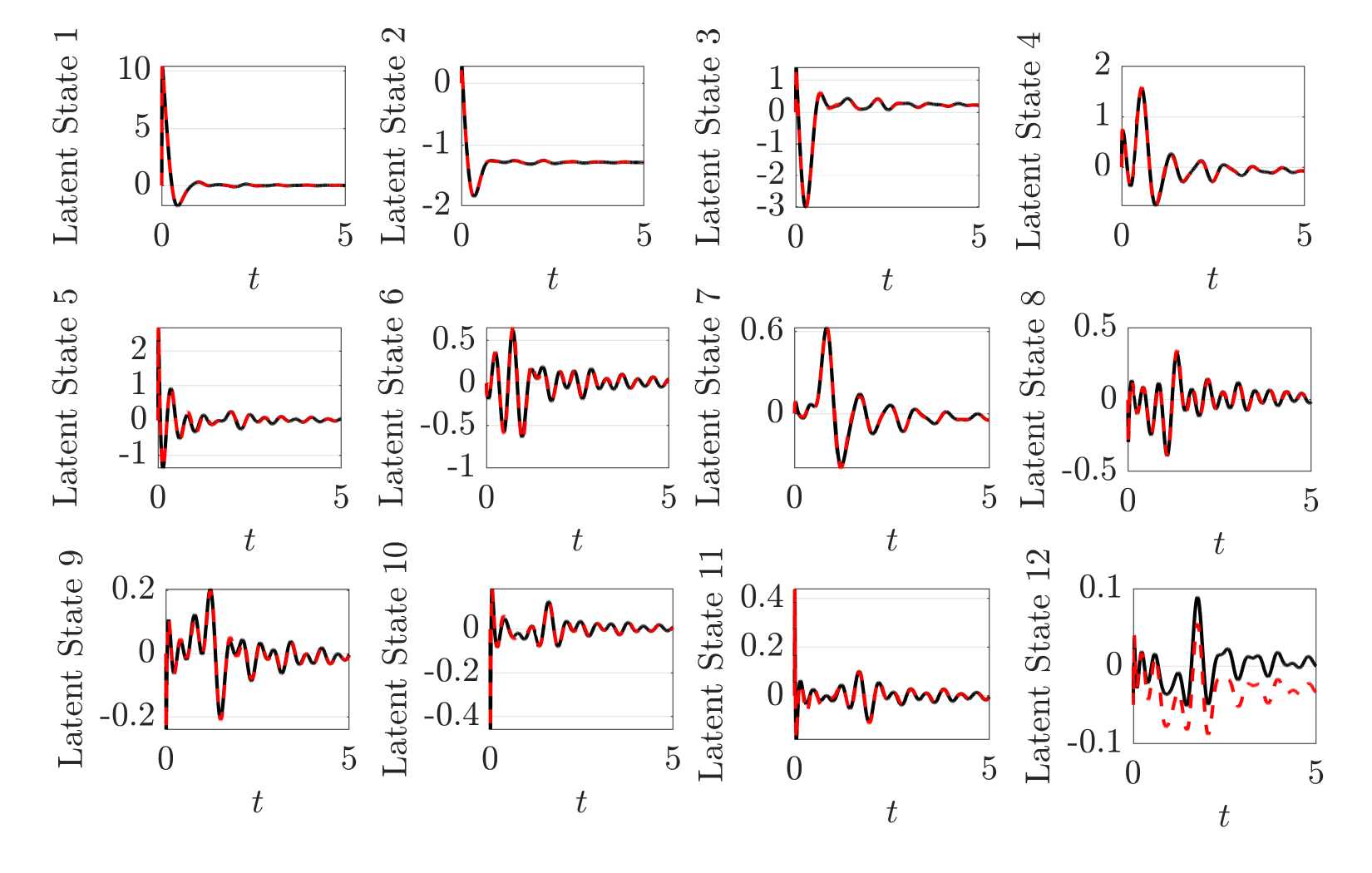}
     \caption{Original latent states and their SINDy predictions for the IEEE 118 Bus systems. The solid lines denote the actual latent states obtained by compressing the true data, and the dashed lines denote the SINDy predictions.}
    \label{fig:Latent118}
\end{figure}

\begin{figure}[H]
  \centering
\includegraphics[width=0.7\linewidth]{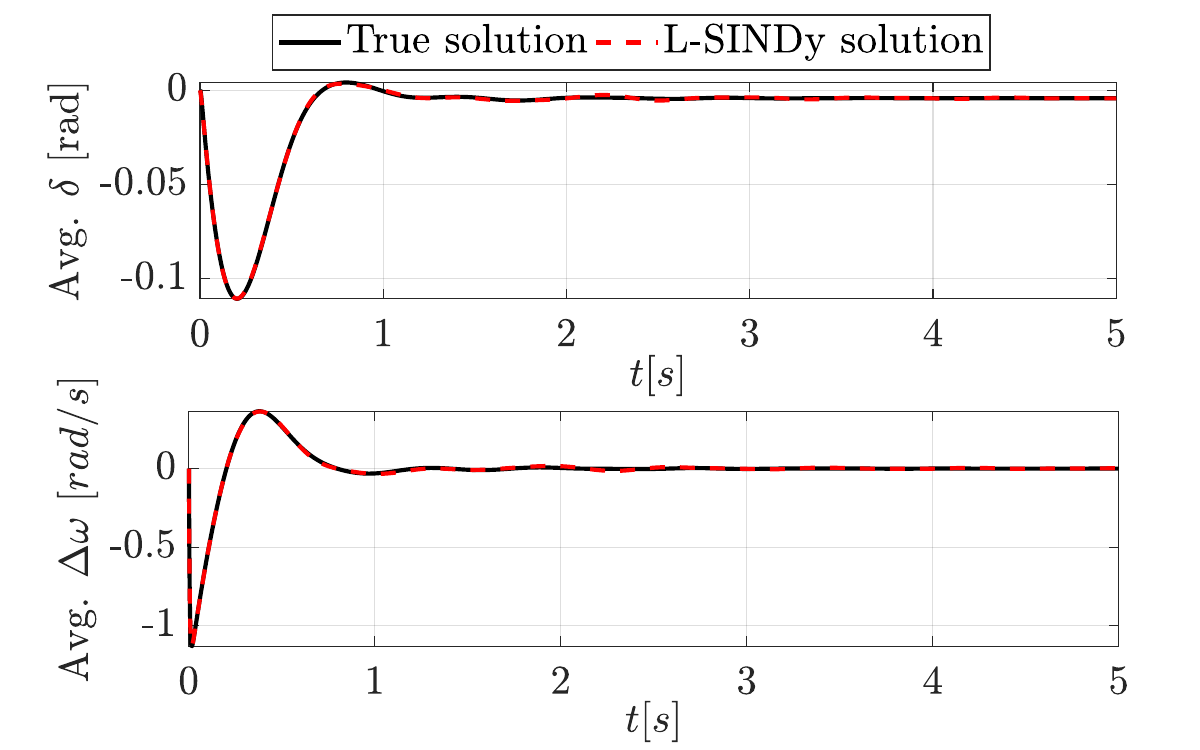}
    \caption{Comparison of true and estimated values of avg. $\delta$ and avg. $\Delta\omega$ in IEEE-118 bus system.}
    \label{fig:pred118}  
\end{figure}

\begin{figure}[H]
    \centering
\includegraphics[width=0.6\linewidth]{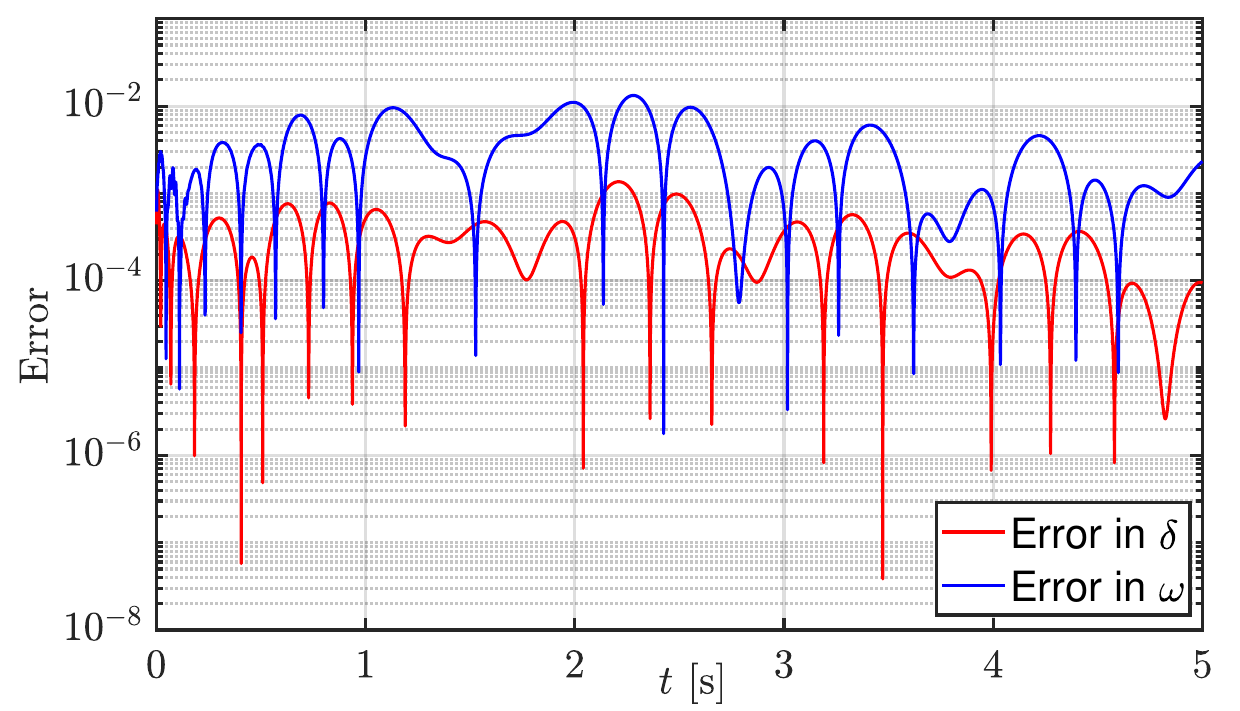}
    \caption{Error plots of the average $\delta$ and avg. $\omega$ error in IEEE-118 bus system.}
    \label{fig:error118}
\end{figure}

The IEEE 118 Bus model is a benchmark model comprising 118 buses and $n_g=54$ generators. To test the L-SINDy method, $m=500$ time snapshots were recorded, consisting of rotor angles $\delta_i$ and corresponding frequencies $\omega_i$ of all the generators. Using the SVD on the data- matrix $\mathbf{X}$, we obtained the singular values. The decay of singular values and the energy captured by each singular value is presented in Figs. \ref{fig:decay118} and Fig. \ref{fig:energy118}, respectively. From the figures, we can observe the slow decay of singular values, a typical characteristic of power system models. Using Algorithm \ref{alg:L-SINDy}, we identified a reduced state model of size $r=12$ by retaining 12 dominant singular values that capture $\approx 99.9\%$ energy in the data. It should be noted that the choice of reduced dimension plays a significant role in the overall quality of the L-SINDy method. In Fig. \ref{fig:Latent118}, we
present the comparison of all the 12 latent states against the latent states obtained by compressing the true data. As can be seen, L-SINDy correctly identifies the latent evolution of these variables. In Fig. \ref{fig:pred118}, we show the time evolution of the average of all the rotor angles $\delta$ and frequency deviations $\Delta \omega$ for the true solution and the L-SINDy solution. As can be seen, L-SINDy accurately predicts the system dynamics, which is also verified from the $\ell_{2}-$norm errors plotted in Fig.~\ref{fig:error118} and reported in Table \ref{tab:comparison}.
\subsection{IEEE-300 Bus system}
\begin{figure}[H]
\centering
\includegraphics[width=0.6\linewidth]{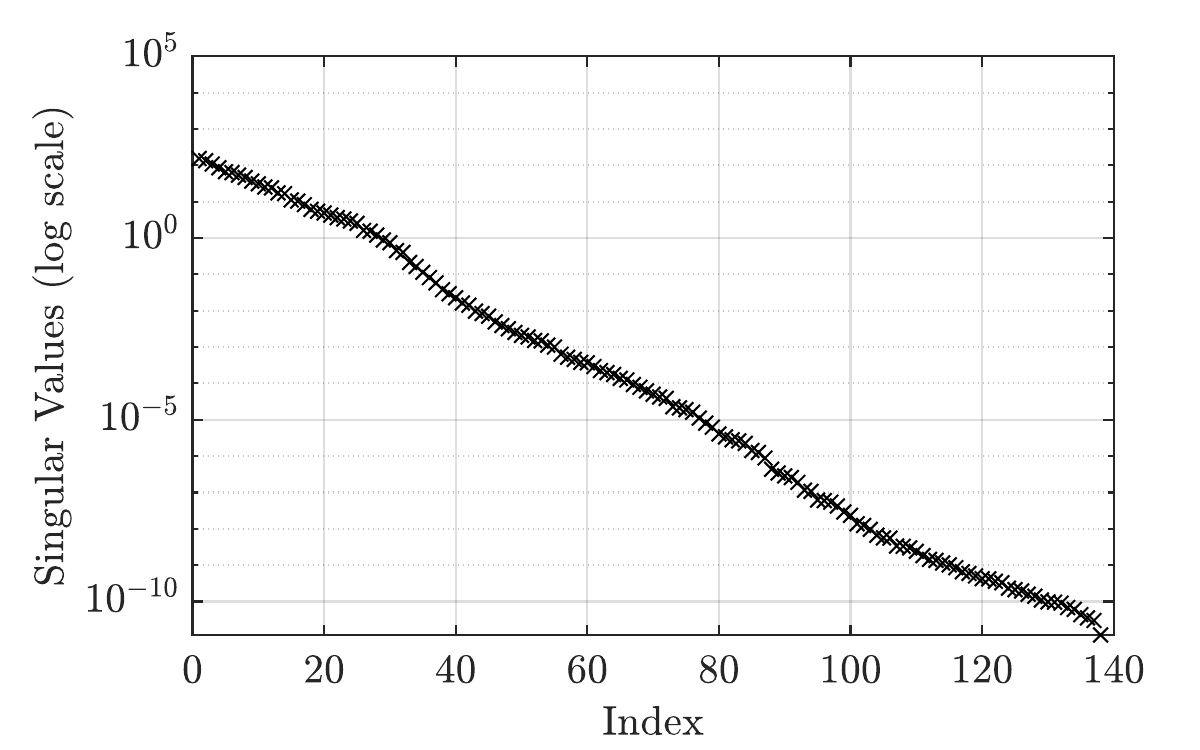}
		\caption{Singular value decay plot of IEEE-300 bus system}
		\label{fig:decay300}
\end{figure}

\begin{figure}[H]
\centering
\includegraphics[width=0.6\linewidth]{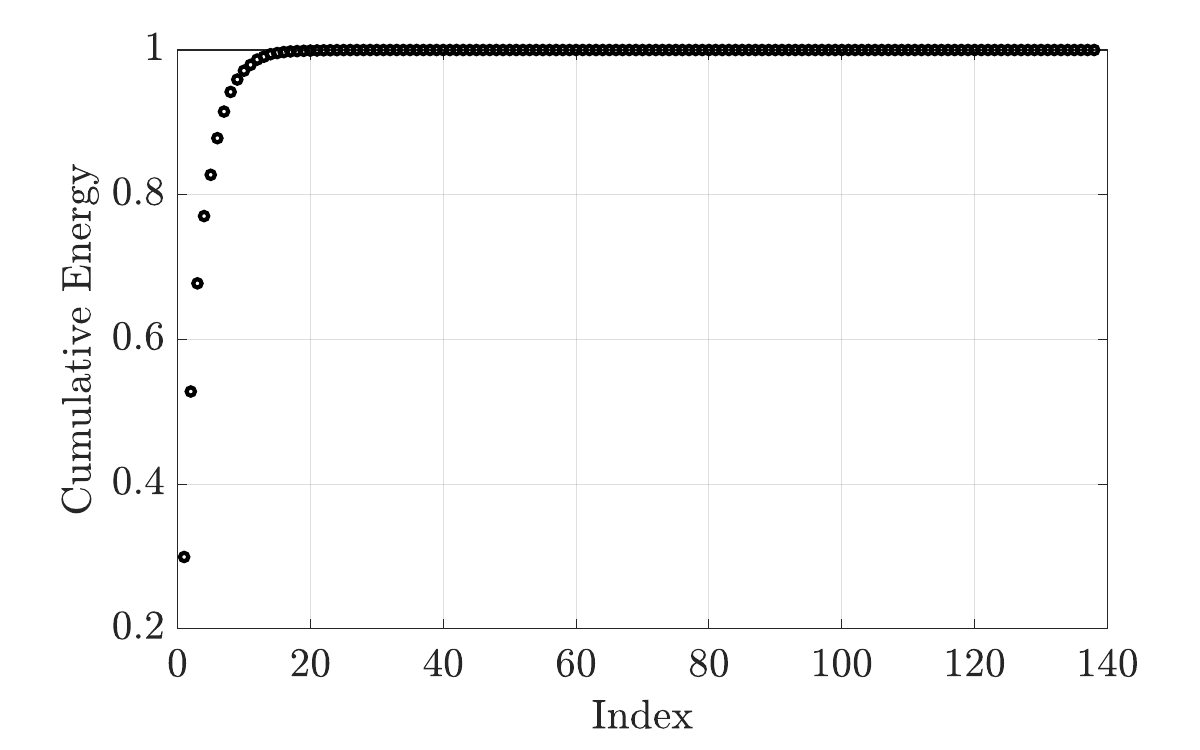}
	\caption{Singular value energy plot of IEEE-300 bus system}
	\label{fig:energy300}	
\end{figure}

\begin{figure}[H]
    \centering
\includegraphics[width=\linewidth]{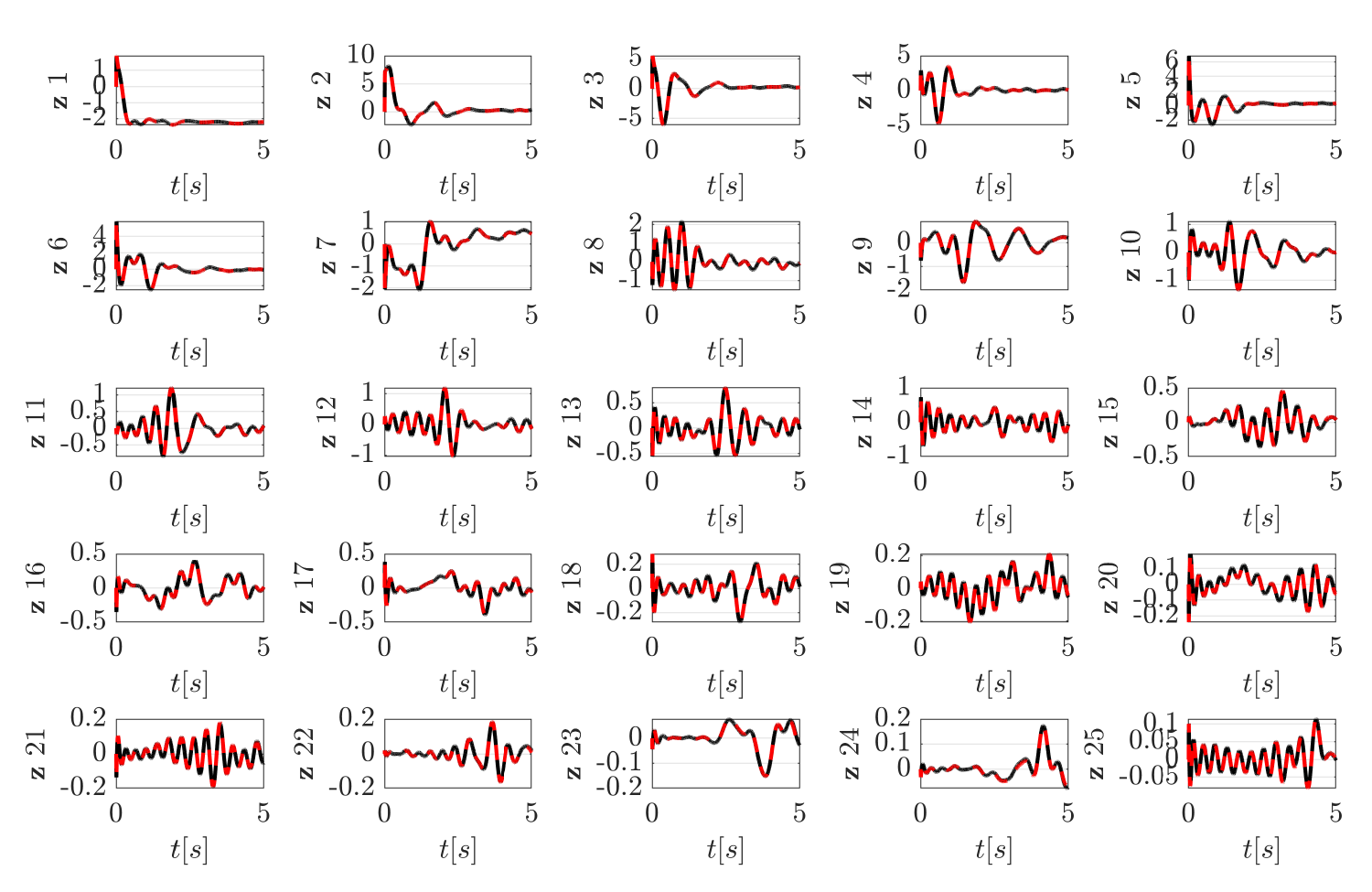}
    \caption{Original latent states and their SINDy predictions for the IEEE 300 Bus systems. The solid lines denote the actual latent states obtained by compressing the true data, and the dashed lines denote the SINDy predictions.}
    \label{fig:300latent}
\end{figure}

\begin{figure}[H]
    \centering
    \includegraphics[width=0.6\linewidth]{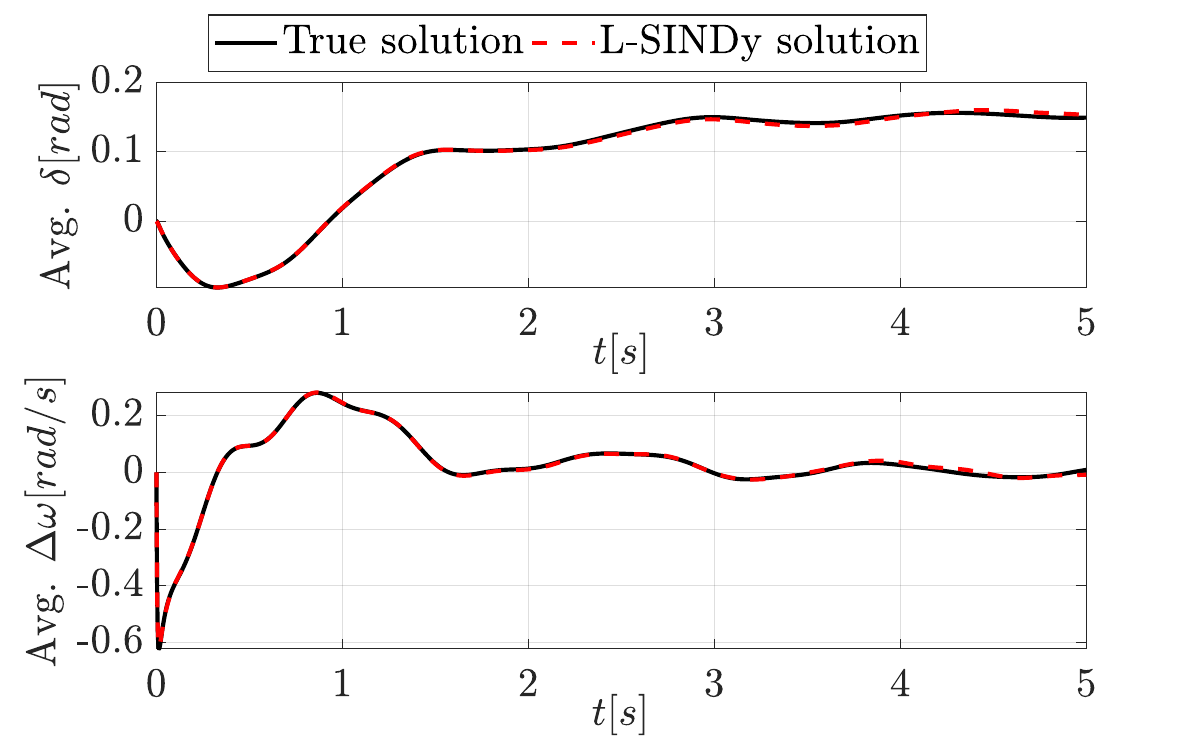}
    \caption{Comparison of true and estimated values of avg. $\delta$ and avg. $\Delta\omega$ in IEEE-300 bus system.}
    \label{fig:300_pred}
\end{figure}

\begin{figure}[H]
    \centering    \includegraphics[width=0.6\linewidth]{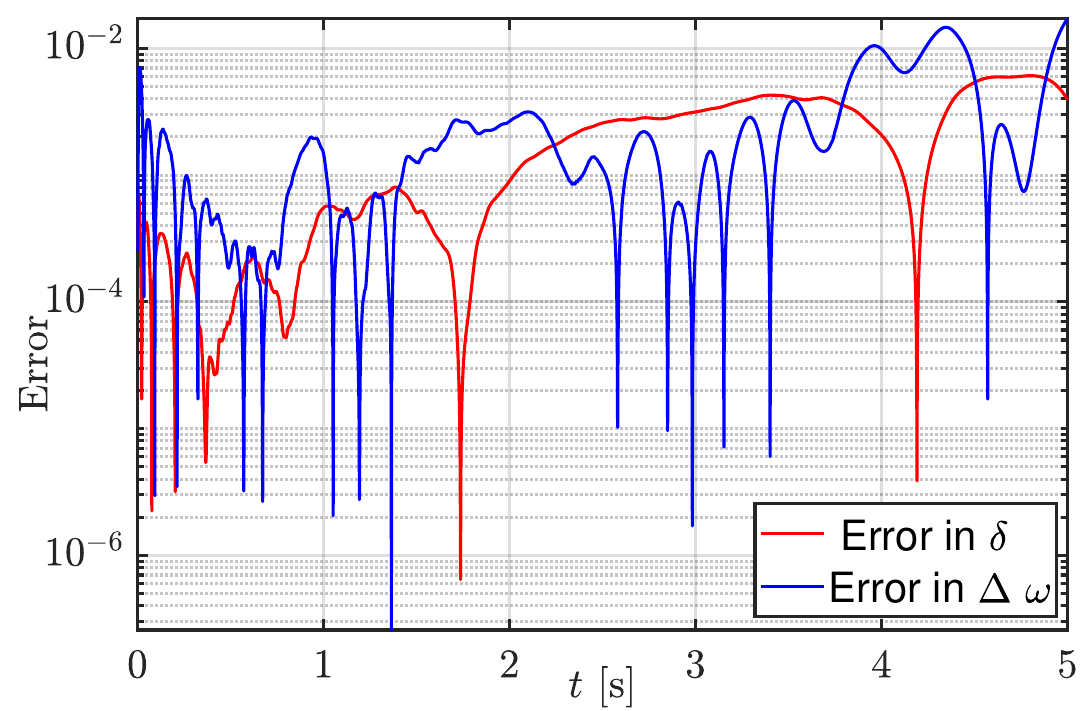}
    \caption{Comparison of avg. $\delta$ and avg. $\omega$ error in IEEE-300 bus system.}
    \label{fig:300error}
\end{figure}

Next, we used the IEEE-300 Bus system. This system represents a considerably more complex test case compared to the previous model. The decay of singular values and the energy captured for this case are shown in  Figs. \ref{fig:decay300} and \ref{fig:energy300} respectively. We can observe that the system exhibits a much slower decay, necessitating a higher size of ROM. For this case, we used $r=25$ POD modes while using Algorithm \ref{alg:L-SINDy}. The plots of latent variables are shown in Fig. \ref{fig:300latent}, and the output response is presented in Fig. \ref{fig:300_pred}. The error plot of the output is shown in Fig. \ref{fig:300error} and quantified in Table \ref{tab:comparison}. The results indicate that the errors remain consistently small and stable over time, demonstrating the superior performance of the proposed L-SINDy method.


\subsection{European High Voltage transmission system}
\begin{figure}[H]
\centering
\includegraphics[width=0.6\linewidth]{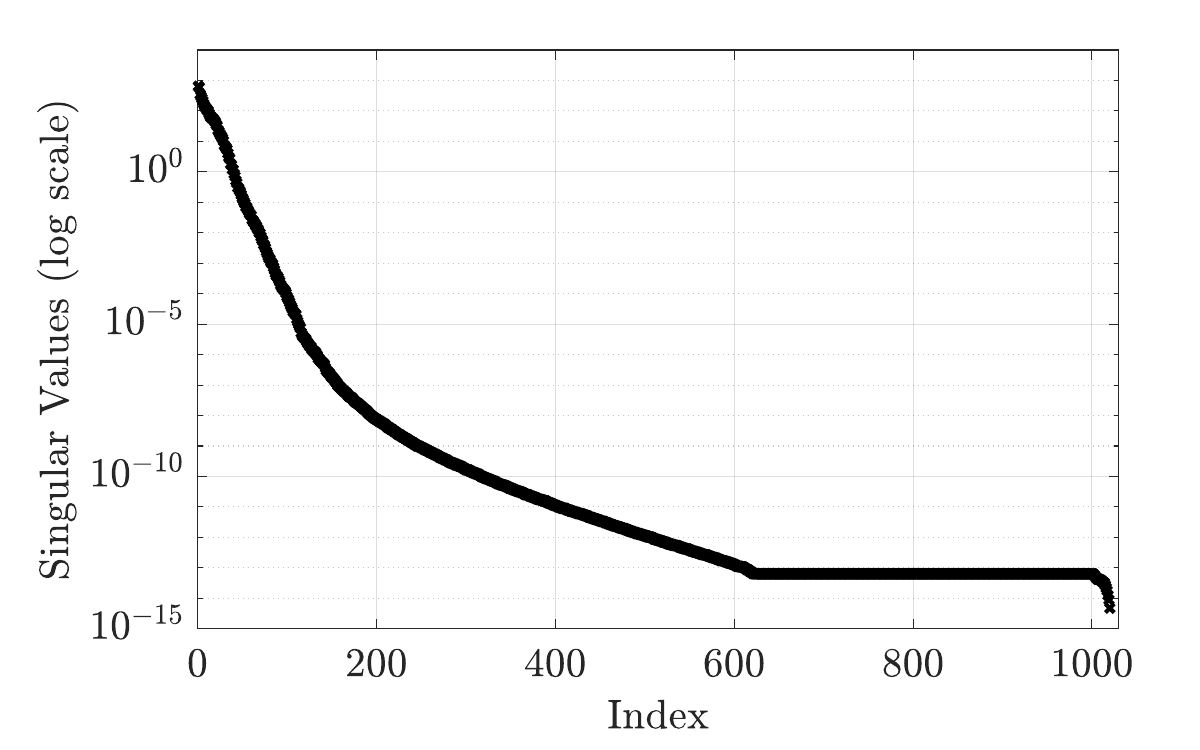}
	\caption{Singular value decay plot of European High Voltage  system}
	\label{fig:decay2869} 
    \end{figure}

\begin{figure}[H]
    \centering
\includegraphics[width=0.6\linewidth]{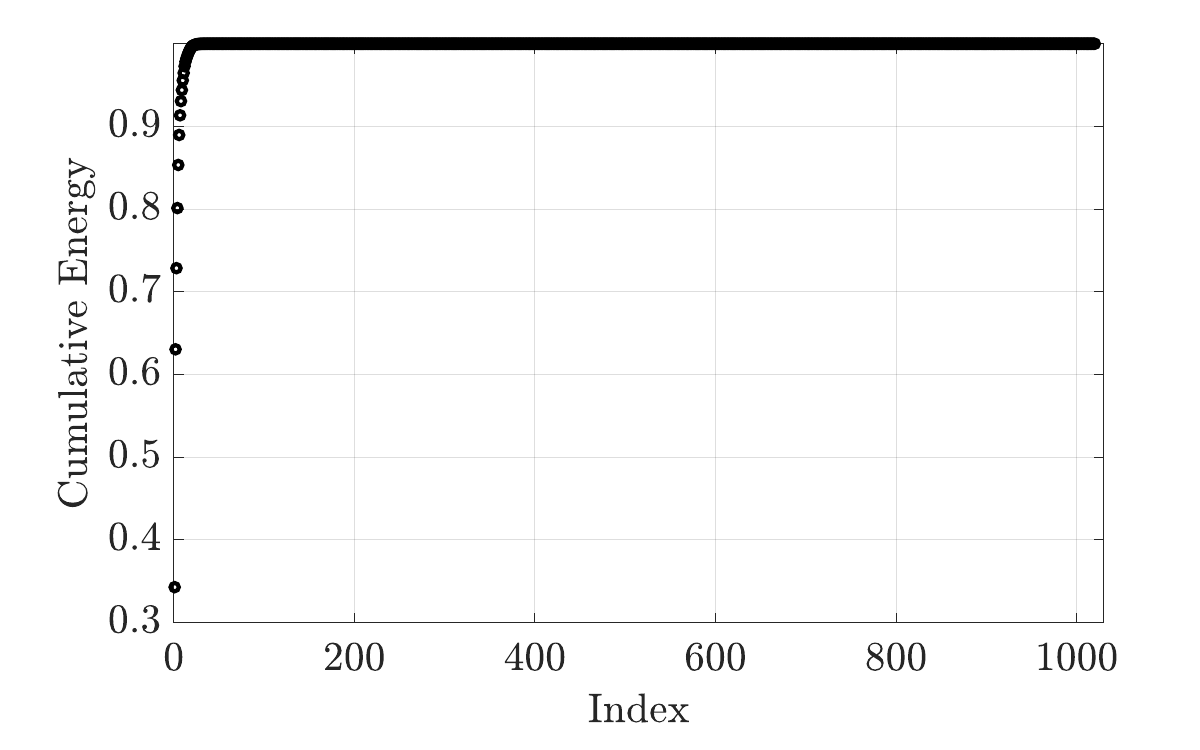}
\caption{Singular value energy plot of European High Voltage  system}
\label{fig:energy2869}		
\end{figure}

\begin{figure}[H]
    \centering
    \includegraphics[width=\linewidth]{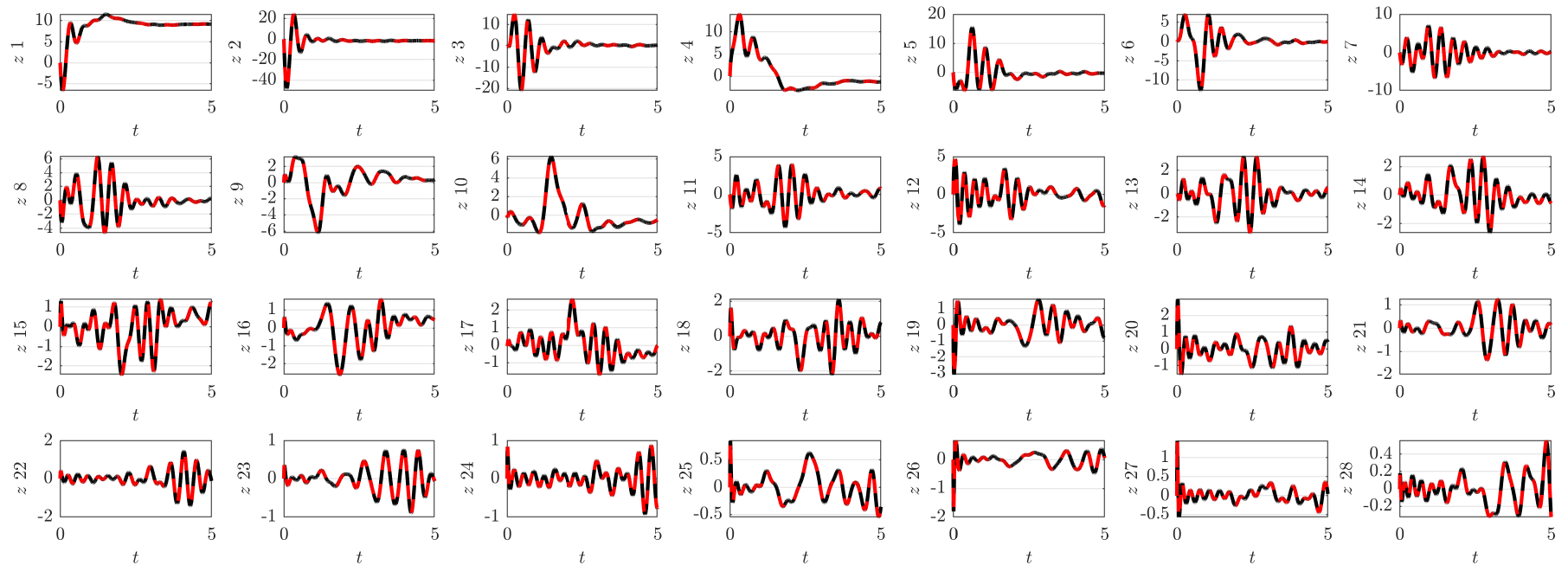}
    \caption{Original latent states and their SINDy predictions for the IEEE 2869 Bus systems. The solid lines denote the actual latent states obtained by compressing the true data, and the dashed lines denote the SINDy predictions.}
    \label{fig:2869_latent}
\end{figure}

\begin{figure}[H]
    \centering   
    \includegraphics[width=0.6\linewidth]{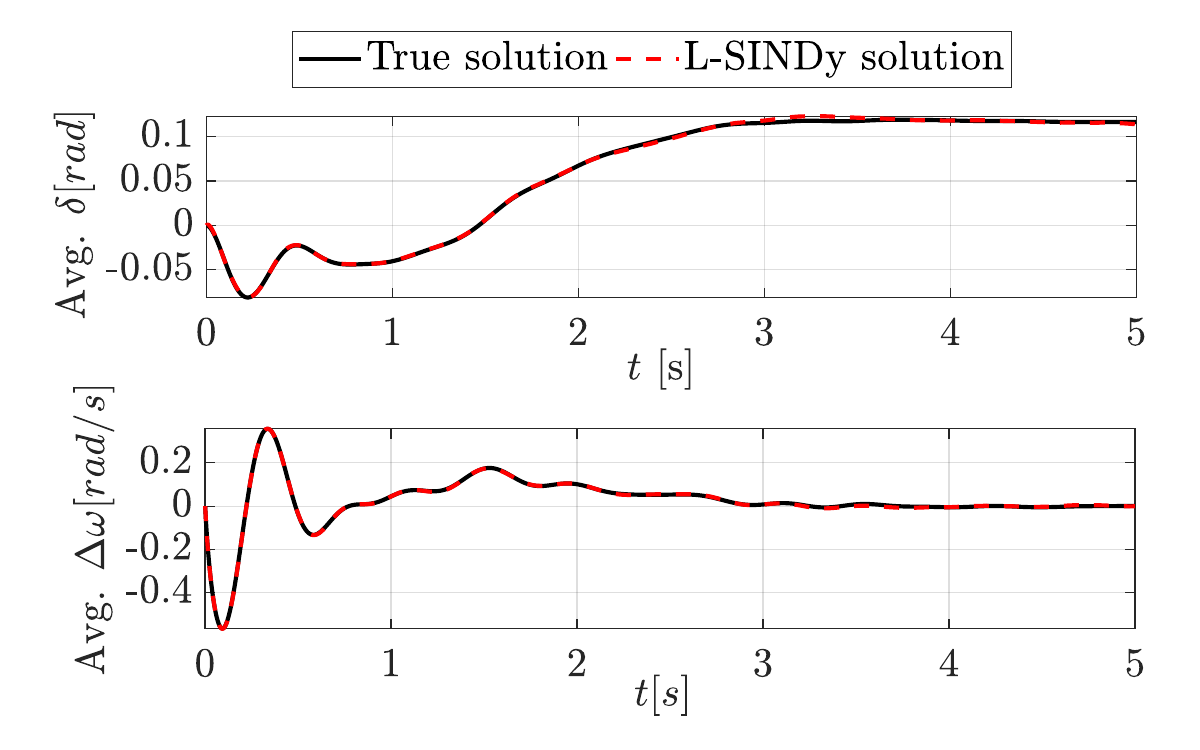}
    \caption{Comparison of true and estimated values of avg. $\delta$ and avg. $\Delta\omega$ in the European High Voltage  system.}
    \label{fig:2869_pred}      
\end{figure}

\begin{figure}[H]
    \centering
    \includegraphics[width=0.6\linewidth]{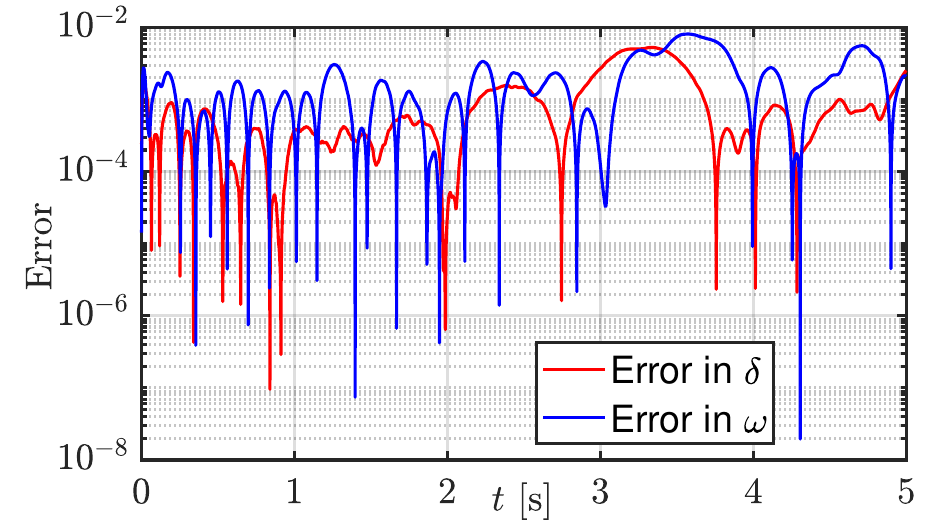}
    \caption{Comparison of avg. $\delta$ and avg. $\omega$ error in the European High Voltage  system.}
    \label{fig:2869_error}
\end{figure}

	Finally, to check the scalability of the proposed L-SINDy method, we used a real-sized power system. The model consists of 2869 buses and 510 generators, resulting in the model dimension of $\approx 2^{10}$. The decay of singular values and the energy captured by POD models is presented in Figs. \ref{fig:decay2869} and \ref{fig:energy2869} respectively. Using Algorithm \ref{alg:L-SINDy}, we determined an appropriate reduced state model of size $r=28$ which captures $\approx 99.9\%$ energy in the system.  Further analysis of the latent-state trajectories, as presented in Fig. \ref{fig:2869_latent}, highlights how the reduced-order model captures the essential coherent patterns within the 2869-bus system. Despite the complexity and size of the network, only a small number of modes were sufficient to describe the coherent generator dynamics. The output plot showing the true and L-SINDy is shown in Fig.~\ref{fig:2869_pred}, and the corresponding error plot is shown in Fig.~\ref{fig:2869_error}. We observe from the figures that L-SINDy can capture the strong nonlinear behavior of the system effectively. This is also verified from the relative $\ell_2$-norm errors reported in Table \ref{tab:comparison}.
    

\begin{table}[htbp]
 \caption{Comparison of relative $\ell_2$ norm errors and CPU time}
    \centering
    	\resizebox{\textwidth}{!}{%
    \begin{tabular}{lcccccc}
    \toprule
        \textbf{Case Study} & {Polyorder}  & $ \lambda $  & $ err_\delta $  & $  err_\omega$ & {FOM-Time $(s)$} & {ROM-Time$(s)$} \\
        \midrule
      \begin{tabular}{l}
           IEEE 118 Bus  \\
            $(n_g =54, n=108, r=12)$ 
      \end{tabular}
         & 4 & 0.001 & $1.6\times 10^{-2}$ & $3\times 10^{-2}$ & 6.87 & 2.31 \\    
       \begin{tabular}{l}
           IEEE 300 Bus  \\
            $(n_g= 69, n=138, r=25)$ 
      \end{tabular}   & 1 & 0.001 & $2.3\times 10^{-2}$  & $3.5\times 10^{-2}$ & 8.42 & 4.12\\
       \begin{tabular}{l}
           IEEE 2869 Bus \\
            $(n_g=510, n=1020, r=28)$ 
      \end{tabular}  & 1 & 0.001 & $1.8\times 10^{-2}$ & $2.3\times 10^{-2}$ & 12.65& 7.54\\
      \bottomrule
    \end{tabular}
   }
    \label{tab:comparison}
\end{table}

\section{Conclusion and Future Work}\label{sec-conclusion}
This paper presents a novel SVD-assisted SINDy framework for estimating dynamic parameters in large-scale low-inertia power grids. By using a reduced-order modeling framework with sparse regression, the technique can identify power grid dynamics accurately and efficiently. This is achieved by solving the sparse regression problem on the reduced POD modes instead of the original state variables. This drastically reduces the size of the variables for the sparse optimization and yields significant computational savings. We demonstrate our findings on the IEEE 118-bus, 300-bus, and 2869-bus test cases. The identified sparse models closely match the true system dynamics, confirming the viability of the proposed L-SINDy method. We provide numerical algorithms and the source codes for researchers and practitioners. \par
The future work will involve testing the proposed scheme for renewable-rich grids, where the dynamics will evolve more rapidly. Also, the L-SINDy method can be made more robust by introducing uncertainty quantification measures to test the viability of the surrogate models during inference. Finally, integrating L-SINDy with online learning from PMU data would enable real-time monitoring of inertia changes in large-scale grids.

\section{Data availability}
The code to reproduce the results in this paper is freely available at GitHub: \url{https://github.com/DanishRaf32/LSINDy}

\bibliographystyle{unsrt}
\bibliography{references.bib}

\end{document}